\documentclass[a4paper,12pt]{article}
\usepackage{amsmath,amsfonts,amssymb,latexsym,ctable}
\usepackage{theorem}
\usepackage{verbatim}
\usepackage{dcolumn}
\usepackage{setspace}
\usepackage{fullpage}
\usepackage{caption}
\usepackage{graphicx}
\usepackage{lineno}
\newtheorem{proposition}{Proposition}

\newtheorem{corollary}{Corollary}

\newtheorem{lemma}{Lemma}
\newtheorem{theorem}{Theorem}
\def\sumn{\underset{i=1}{\overset{N}{\sum}}}
\def\sumkzero{\underset{k=0}{\overset{\infty}{\sum}}}
\def\sumkone{\underset{k=1}{\overset{\infty}{\sum}}}

\def\un {\hbox{\bf 1}}
\def\ex {\mathbb{E}}
\def\var{\mathbb{V}}

\begin{document}
\parindent 0pt  
\renewcommand{\baselinestretch}{1.5}

\title{Long-memory process and aggregation of AR(1) stochastic processes: A new characterization}
\date{} 
\author{Bernard Candelpergher$^{a}$, Michel Miniconi$^{b}$ and Florian Pelgrin$^{c}$}
\maketitle

\begin{abstract}
\noindent  Contemporaneous aggregation of individual AR(1) random processes might lead to different properties of the limit aggregated time series, in particular, long memory (Granger, 1980). We provide a new characterization of the series of autoregressive coefficients, which is defined from the Wold representation of the limit of the aggregate stochastic process, in the presence of long-memory features. 
Especially the infinite autoregressive stochastic process defined by the almost sure representation of the aggregate process has a unit root in the presence of the long-memory property. Finally we discuss some examples using some well-known probability density functions of the autoregressive random parameter in the aggregation literature.
\medskip

\noindent JEL\ Classification Code: C2, C13.
\medskip

\end{abstract}

\vspace*{2cm}

\vspace*{0.4cm} \noindent$\mbox{}$$^{a}${\footnotesize University of Nice Sophia-Antipolis, Laboratoire Jean-Alexandre Dieudonn\'{e}. E-mail: bernard.candelpergher@unice.fr.}

\medskip\noindent$\mbox{}$$^{b}${\footnotesize University of Nice Sophia-Antipolis, Laboratoire Jean-Alexandre Dieudonn\'{e}. E-mail: michel.miniconi@unice.fr.}

\medskip\noindent$\mbox{}$$^{c}${\footnotesize EDHEC Business School, E-mail: florian.pelgrin@edhec.edu.}

\vspace*{0.2cm} {\footnotesize \medskip\noindent This paper has benefited from very useful comments from and discussions with St\'{e}phane Gregoir. The usual disclaimer applies. }

\newpage
\normalsize

\section{Introduction}

\noindent \noindent Aggregation is a critical and widely acknowledged issue in the empirical and theoretical literature in economics and other fields. Especially, since the contributions of Granger and Morris (1976) and Granger (1980), it is well-known that the contemporaneous aggregation of individual random AR(1) stochastic processes might lead to long memory models.\footnote{In the sequel, a stochastic process is said to have a long memory property if its autocovariance function is not summable (Beran, 1994; Beran et al., 2013).} Notably, Granger (1980) considers the case of a Beta-distribution for the random autoregressive parameter and thus points out that the long memory  property depends on the behavior of the density of the random autoregressive parameter near unity and that the common and idiosyncratic components might exhibit a different degree of long memory. In particular, Zaffaroni (2004) generalizes  these results by studying the limit of the aggregate process with a quite flexible (semi-parametric) assumption regarding the behavior near unity of the probability density function of the random autoregressive coefficient and makes clear the asymptotics of both the common and idiosyncratic parts of the aggregated process. Among others results, Zaffaroni (2004) shows formally that the more concentrated is the distribution of the random autoregressive coefficient near the unit, the stronger is the long-memory property of the limit aggregated process.\footnote{For a recent survey, see Leipus et al. (2013).}. Following these contributions, we study the aggregation of heterogenous individual AR(1) stochastic processes that leads to the long-memory property. In contrast to the literature, we focus on the infinite autoregressive representation of the limit of the aggregate process, and especially the sum of the autoregressive coefficients, rather than the usual infinite moving representation.
\\ \\
\noindent Indeed the use of the (equivalent) infinite autoregressive representation of the limit of the aggregated process might be meaningful in various contexts of the aggregation (panel) literature. As pointed out by Lewbel (1994), the autoregressive representation is insightful among others  when estimating the (macro) aggregate dynamics with unobservable individual series in the presence of unobserved parameter heterogeneity, when identifying and estimating certain distributional features of the micro parameters from aggregate relations (disaggregation problem), or when deriving the aggregate dynamics with unknown common error terms.\footnote{See also Pesaran (2003), Pesaran and Chudik (2014), Jondeau and Pelgrin (2014a,b).} Obviously, if the limit of the aggregate process is represented almost surely by a short memory process, the infinite autoregressive representation and especially the determination of its persistence (e.g., through the sum of the autoregressive coefficients) easily obtains using standard results of time series analysis (Brockwell and Davis, 2002). In contrast, if the limit of the aggregate process is represented almost surely by a long memory process, the convergence and thus the determination of the series of autoregressive coefficients is challenging since the series of moving average coefficients is no longer absolutely summable. Such a characterization might be necessary for instance in the case of the estimation of the aggregate dynamics with unobservable individual series in which finite parameter approximation for the infinite lag distribution is required---the autoregressive representation of the limit aggregated process displays an infinite parametrization whereas the econometrician has only finitely many observations and thus finite-parameter approximations might lead to near observational aggregate processes with different statistical properties as long as the sum of the autoregressive coefficients is not correctly identified. In this respect, our paper tackles this issue and proposes a new characterization of the autoregressive representation of limit long-memory aggregate processes that results from individual random AR(1) stochastic processes.
\\ \\
\noindent Notably aggregation of individual random AR(1) processes is analyzed in this paper under the assumptions that some common factors exist, the (positive) random autoregressive coefficient takes values in $[0,1)$ with probability distribution $\mu$, and all of the noncentral moments exist. In particular, the series of all of the noncentral moments might be either absolutely summable or divergent and the limit of the aggregated stochastic process satisfies almost surely both an infinite moving average and autoregressive representations. Within this framework, 
we show that the sum of the autoregressive coefficients of the stochastic process defined from the limit in $\mathcal{L}^2$ of the aggregate process, say $\underset{k \geq 1}{\sum}a_k$, might equal one (respectively, less than one) when the limit aggregate process has the long memory property (respectively, short memory property). Say differently, the divergence of the series of noncentral moments is fully equivalent to the presence of a unit root in the representation of the stochastic process defined by the limit of the aggregate process.  In so doing, we consider the mapping between the moving average coefficients and the autoregressive coefficients of the limit aggregate process and make use of complex analysis when the series of the moving average coefficients is divergent. Complex analysis is called for at least two arguments. On the one hand, the existence of the limit 
$\lim_{r\rightarrow 1^{-}}{a(r)}$ of the 
(ordinary) generating function of the sequence of the autoregressive coefficients does not insure its equality to $\sum_{k\geq 1} a_k$ or even the convergence of this sum. On the other hand, we cannot apply a standard Hardy-Littlewood Tauberian theorem since, to the best of our knowledge, there is no general proof of the positiveness of the autoregressive coefficients for all probability distributions $\mu$ defined on $[0,1)$. Interestingly such a complex analysis might be used to study the behavior on the unit circle of an infinite parameterized stochastic process, which is defined (almost surely) in $\mathcal{L}^2$ but does not belong to $\mathcal{L}^1$. 
\\ \\
\noindent The rest of the paper is organized as follows. In Section 2, we discuss the main assumptions regarding the individual random AR(1) stochastic processes and then we derive the limit aggregate process. In Section 3, we provide the two main results of our paper, namely the determination of the sum of the autoregressive parameters and the characterization of the stochastic process defined by the limit in $\mathcal{L}^2$ of the aggregated process. In Section 4, we assume that the distribution of the random autoregressive coefficient is subsequently a Beta distribution (of the first kind), a uniform distribution and a random polynomial density function. In so doing, we apply the results of Section 3 and generalizes some results in the literature. Proofs are gathered in the Appendix.

\section{Aggregating individual AR(1) processes}

In this section, we first discuss the assumptions regarding the individual AR(1) processes. Then we derive the limit aggregate process.

\subsection{Individual dynamics}

\noindent Consider the individual random AR(1) model for $i=1,\cdots,N$ and $t \in \mathbb{Z}$:\footnote{Alternatively one might assume that $t \in \mathbb{N}$ and consider an asymptotically stationary limit aggregate process (Gon\c{c}alves and Gourieroux, 1994).}
\begin{equation}
x_{i,t}=\varphi_{i} \ x_{i,t-1}+v_{i,t}, \label{xi,t}
\end{equation}%
where $\varphi_i$ denotes an individual-specific (random) parameter drawn from a fixed random variable $\varphi$ and $v_{i,t}$ is an error term that can decomposed into a common
component, $\epsilon_{t}$, and an idiosyncratic (individual-specific) component, $\eta_{i,t}$:\footnote{Such dynamics have been used in economics (among others) to represent consumption expenditures across households (Lewbel, 1994), consumer price inflation across subindices (Altissimo et al., 2009), real exchange rates across sectors (Imbs et al., 2005), or real marginal cost across industries (Imbs et al., 2011).}
\begin{eqnarray}
v_{i,t} &=&\epsilon_{t} + \eta_{i,t}\label{vi,t}.
\label{ui,t}
\end{eqnarray}%
\noindent The macro variable results from the aggregation of micro-units, with the use of time-invariant nonrandom weights $W_N = \left(w_1,\cdots,w_N\right)'$, with $\sumn w_i = 1$, so that the aggregate variable is defined as the weighted average of the micro-units  $X_{N,t} = \sumn w_i x_{i,t}$. The following assumptions hold:
\newline \newline
\noindent \textbf{Assumption 1:} $\varphi$ a fixed random variable with probability distribution 
$\mu$ the support of which is in $[0,1]$.
\newline
\noindent \textbf{Assumption 2:} The moments $u_k=\ex(\varphi^{k})$ exist for all integer $k\geq 1$.
\newline
\noindent \textbf{Assumption 3:}  $\epsilon_t$ and $\eta_{i,t}$ are white noise processes with means of zero and variance of $\sigma_{\epsilon}^{2}$ and $\sigma_{\eta}^{2}$, respectively; $\epsilon_t$ and $\eta_{i,t}$ are mutually orthogonal at any lag and lead.
\newline
\noindent \textbf{Assumption 4:} Realizations of $\varphi$ are independent of those of $\epsilon_t$.
\newline
\noindent \textbf{Assumption 5:} As $N \rightarrow \infty$, $\left\|W_N\right\| = O\left(N^{-1/2}\right)$ and $w_i/\left\|W_N\right\| = O\left(N^{-1/2}\right)$ for all $i \in \mathbb{N}$.
\newline \newline
\noindent Before discussing our main assumptions, we introduce some notations. Let $\mathcal{H}_x$ denote the Hilbert space generated by all of the random variables that compose the panel $(x_{i,t})_{i=1,\cdots,N, t \in \mathbb{Z}}$, $\mathcal{H}_{x_i}$ the Hilbert space generated by the stochastic process $(x_{i,t})_{t \in \mathbb{Z}}$, $\mathcal{H}_{x,t}$ and $\mathcal{H}_{x_i,t}$ the corresponding subspaces of $\mathcal{H}_{x}$ and $\mathcal{H}_{x_i}$ up to time $t$. Assumption 1 outlines that $\varphi$ is a random variable with distribution on the support $[0,1)$ such that $\varphi \in \underset{t}{\bigcap}\underset{i}{\oplus} \mathcal{H}_{x_i,t}$.\footnote{Note that the set $(-1,0)$ has been excluded from the support of $\varphi$ for simplicity's sake.} This assumption is consistent with many parametric specifications of the cross-sectional distribution of $\varphi$. We only rule out situations in which some individual processes are not almost surely (asymptotically) stationary, i.e. $\mathbb{P}\left(|\varphi| \geq 1\right) > 0$ (see further). In particular, Assumption 1 includes the Beta distribution, $B(p,q)$. In this case, the representation of the (limit) aggregate process does depend on the mass distribution of the Beta distribution around unity: the smaller is $q$, the larger is the mass of the distribution around unity (Zaffaroni, 2004). In contrast, imposing the condition $0 \leq \varphi \leq c < 1$ for some constant $c$ will guarantee that there are no individual unit root parameters that would dominate at the aggregate level (Zaffaroni, 2004) and that the limit aggregate process (as $N \rightarrow \infty$) displays short memory with an exponentially decaying autocorrelation function. Assumption 2 insures that noncentral moments $u_k=\mathbb{E}(\varphi^{k})$ of any (nondegenerate) random variable $\varphi$, defined on $[0,1)$, satisfy: $1>u_{1}\geq \cdots \geq u_{k} \geq 0, \forall k\geq 1,$ and $u_{k} \rightarrow 0$ as $k\rightarrow \infty$. 
\\ \\
\noindent Assumption 3 defines the statistical properties of the two components of the error term, the common shock (factor) and the idiosyncratic shock. Several points are  worth commenting. First, without loss of generality, one might assume that the stochastic process $(v_{i,t})_{i=1,\cdots,N, t\in \mathbb{Z}}$ is weakly linearly exchangeable (Aldous, 1980), i.e. $(v_{i,t})_{i=1,\cdots,N, t \in \mathbb{Z}}$ is a purely non deterministic (regular) covariance stationary process and the covariance structure is invariant by translation with respect to the time index and invariant by permutation with respect to the unit index. Second, the stochastic process $(\epsilon_{t})_{t \in \mathbb{Z}}$ is assumed to be known and given.\footnote{Without altering our main results but at the expense of further assumptions, (i) the common error term at time $t$ might be multiplied by the realization of a (scaling) random variable $\kappa$ to introduce some form of heteroscedasticity, (ii) the stochastic process $(\epsilon_{t})_{t \in \mathbb{Z}}$ might be assumed to be unknown (using a two-step procedure), and (iii) multiple independent common error terms might be introduced in the specification (Zaffaroni, 2004).} Given that $\mathcal{H}_{x} = \underset{i}{\oplus}\mathcal{H}_{x_i}$ and thus $\mathcal{H}_{x,t} = \underset{i}{\oplus}\mathcal{H}_{x_i,t}$, one has $\epsilon_t \in \mathcal{H}_{x,t} \bigcap (\mathcal{H}_{x,t-1})_{\perp}$. Third, taking that $\mathcal{H}_{x_i} = \mathcal{H}_{\epsilon} \oplus \mathcal{H}_{x_i - \epsilon}$, the idiosyncratic stochastic process is such that $\eta_{i,t} \in \mathcal{H}_{x_i - \epsilon,t} \bigcap (\mathcal{H}_{x_i - \epsilon,t-1})_{\perp}$. Fourth, the assumption that the variance of the idiosyncratic shocks is the same across individuals might be too stringent in practical examples (e.g., in economics or in finance) and can be generalized by assuming heteroscedastic idiosyncratic error terms.
Eq. (\ref{vi,t}) together with Assumption 3 provide a parsimonious form of (statistical) cross-sectional dependence, which is common in the aggregation literature (Forni and Lippi, 2001; Zaffaroni, 2004). Obviously one may introduce dependence on the idiosyncratic part through cross-correlation and/or spatial correlation (Chudik, Pesaran and Tosetti, 2011; Pesaran and Tosetti, 2009). Such correlation has no effect on our results as long as it is sufficiently weak (e.g., local correlation) so that the error structure belongs to the class of approximate factor structure models (Chamberlain and Rothschild, 1983).\footnote{The condition is that the maximum eigenvalue of the covariance matrix of $\eta_t = (\eta_{1,t},\cdots,\eta_{N,t})'$ remains bounded as the cross-section dimension increases.}.
Assumption 4 avoids any correlation between the error terms and $\varphi$. Assumption 5 is a granularity condition, which insures that the weights used to define the aggregate process are not dominated by a few of the cross-sectional units (Gabaix, 2011; Pesaran and Chudik, 2014).\footnote{Our results extend to the case of (time-varying) stochastic weights. Such an extension requires at least that the weights be distributed independently from the stochastic process defining the random variable.}

\subsection{Aggregate Dynamics}

\noindent 
The empirical cross-sectional moments of $\varphi$ are $\tilde{\ex}_N\left(\varphi^k\right)=\sumn w_i\varphi^k_i$, $\forall k \geq 1$. For sake of simplicity and without loss of generalization with respect to Assumption 4, we assume that $w_i = 1/N$ for all $i$. Consequently, as $N \rightarrow \infty$,  $\tilde{\ex}_N\left(\varphi^k\right) \overset{a.s.}{\rightarrow} u_k$.
\newline \newline
\noindent Using Eqs. (\ref{xi,t})--(\ref{vi,t}), the exact aggregate dynamics can be written as:\footnote{Put differently, it is an ARMA(N,N-1) in the absence of common roots in the individual processes (Granger and Morris, 1976).}
\begin{eqnarray}
\prod_{j=1}^N\left(1-\varphi_jL\right)X_{N,t} &=& \frac{1}{N}\sum_{i=1}^N \prod_{j \neq i}\left(1-\varphi_j L\right)v_{i,t} \label{exact_1}
\end{eqnarray}
or equivalently
\begin{eqnarray}
X_{N,t} &=& \frac{1}{N}\sum_{i=1}^N \left(1 - \varphi_i L\right)^{-1}\epsilon_{t} + \frac{1}{N}\sum_{i=1}^N \left(1 - \varphi_i L\right)^{-1}\eta_{i,t} \label{exact_2}
\end{eqnarray}
where $L$ is the lag operator ($z_{t-1} = L z_t$). 
Taking Eq. (\ref{exact_2}), we can characterize the asymptotic behavior of both the idiosyncratic component and the common component. This is done in the following proposition. Results are known but are reported here for sake of completeness.

\begin{proposition}
\label{Proposition Aggr Dyn} Suppose that Assumptions 1--5 hold. Given the disaggregate model defined in Eqs. (\ref{xi,t})--(\ref{vi,t}), the limit in $\mathcal{L}^2$ of the aggregated process as $N \rightarrow \infty$ satisfies (almost surely) the two equivalent representations:
\begin{eqnarray}
X_{t} &=&\sum_{k=0}^{\infty }u_{k}\ \epsilon_{t-k} \quad\quad\quad\quad\quad \ \text{(MA form)},  \label{Agg_MA} \\
X_{t} &=&\sum_{k=1}^{\infty }a_{k} X_{t-k}+\epsilon_{t} \quad\quad\quad \text{(AR form)},  \label{AR1_Lewbel}
\end{eqnarray}
where $X_{N,t} \overset{\mathcal{L}^2}{\rightarrow} X_t$ and $\tilde{\ex}_N\left(\varphi^k\right) \overset{a.s.}{\rightarrow} u_k = \ex\left(\varphi^k\right)$ as $N \rightarrow \infty$.
The sequence $\left\{a_k, k \geq 1\right\}$ where $a_k = \ex\left[A_k\right]$ satisfies the recurrence relation :
\begin{eqnarray}
A_1=\varphi\ ,\ A_{k+1} = \left(A_k - a_k\right)\varphi \label{recurrence}
\end{eqnarray}
\end{proposition}

\noindent Proof: Gon\c{c}alves and Gouri\'{e}roux (1988), Lewbel (1994). 
\newline \newline
Several points are worth commenting. First, as shown by Robinson (1978), with the exception of a degenerate distribution for $\varphi$ (e.g., Dirac distribution), the dynamics of the limit aggregate process is richer than the individual dynamics because of the nonergodicity of the individual random AR(1) process. 
Second, using the infinite moving average representation and the positiveness of the moments, the (limit) aggregate process displays short memory if $\sum_{k = 0}^\infty u_k < + \infty$ whereas it has a long-memory property if  $\sum_{k = 0}^\infty u_k = + \infty$.
It is worth noting that (i) the (limit) aggregate process is in $\mathcal{L}^2$ but not necessarily in $\mathcal{L}^1$, and (ii) the sum of the autoregressive coefficients is (absolutely) convergent if the aggregate process has short memory (since the spectral representation is unique).
Third, the central moments of the cross-sectional distribution of $\varphi$ can be easily obtained from the infinite autoregressive representation of the aggregate process. This is useful when considering the standard disaggregation problem in statistics. For instance, the first four cross-sectional moments are $\ex\left[\varphi\right] = a_1$, $\var\left[\varphi\right] = a_2$, $S\left[\varphi\right] =\left(a_{3}-a_{1}a_{2}\right) /\left(a_{2}\right)^{3/2}$, and $K\left[\varphi\right] =\left(a_{4}-2a_{1}a_{3}+a_{1}^{2}a_{2}+a_{2}^{2}\right)/\left(a_{2}\right) ^{2}$.
Fourth, Equation (\ref{AR1_Lewbel}) shows that aggregation leads to an infinite autoregressive model for $X_t$ (see Robinson, 1978, Lewbel, 1994). Notably, using Eq. (\ref{recurrence}), the autoregressive parameters $a_{k}$ are nonlinear transformations of the noncentral moments of $\varphi$ and satisfy the following non-homogenous difference equations  (for $k \geq 2$):
\begin{eqnarray*}
a_{k+1} & \equiv & \ex\left[A_{k+1}\right]=u_{k+1} -\sum_{r=1}^{k}a_{r}u_{k-r+1} \\
a_1 &=& \ex\left[\varphi\right].
\end{eqnarray*}
Fifth, the persistence of the stochastic process defined by the limit of the aggregated process can be defined as the sum of the autoregressive coefficients, which is denoted $a(1) = \sum_{k=1}^{\infty} a_k$. More specifically, if the limit of the aggregate process belongs to $\mathcal{L}^1$ and thus $(X_t)_{t \in \mathbb{Z}}$ is a short-memory stochastic process, then it is straightforward to show that
\begin{eqnarray*}
a(1) = \left(\ex\left[\frac{1}{1-\varphi}\right]\right)^{-1}\ex\left[\frac{\varphi}{1-\varphi}\right].
\end{eqnarray*}
or equivalently
\begin{eqnarray*}
\left(1 - a(1)\right)^{-1} = \ex\left[\frac{1}{1-\varphi}\right].
\end{eqnarray*}
In particular, the limit of the aggregate stochastic process has short memory if and only if $\ex\left[\frac{1}{1 - \varphi}\right] < \infty$ (Robinson, 1978; Gon\c{c}alves and Gourieroux, 1994). In the spirit of Zaffaroni (2004), a sufficient condition on the probability density function $h$ of $\varphi$ for the short memory property is that there exists $\alpha \in (0,1)$ and a constant $C$ such that $\underset{x \rightarrow 1^-}{\lim} \frac{h(x)}{(1-x)^{1-\alpha}} = C$. Sixth, we can define the generating function of the $a_k$ terms as follows. Taking the recurrence relation, one can write formally
\begin{eqnarray*}
\sum_{k \geq 1} A_k z^k = \frac{z\varphi}{1 - z\varphi}\left(1 - \sum_{k \geq 1}a_k z^k\right)
\end{eqnarray*}
and thus
$$
\sum_{k\geq 1}a_{k}z^{k}
=\sum_{k\geq 1}u_{k}z^{k}-\sum_{k\geq 1}a_{k}z^{k}\sum_{k\geq 1}u_{k}z^{k}.
$$
Let $m$ denote
$$
m(z)=\sum_{k\geq 1}u_{k}z^{k},
$$
we obtain the formal generating function of the $a_k$ terms
\begin{eqnarray}
a(z) \equiv \sum_{k\geq 1}a_{k}z^{k}
=\frac{m(z)}{1+m(z)}, \label{generating_function}
\end{eqnarray}
and thus the mapping between the generating function of the infinite autoregressive lag polynomial and the one of the infinite moving average representation. As explained in Section 3, Eq. (\ref{generating_function}) is fundamental since our proof is intimately related to the function $\displaystyle\frac{m}{1+m}$.

\section{Aggregate long-memory process and the sum of the autoregressive parameters}

\noindent In this section, we show that the sum of the autoregressive coefficients equals one in the presence of a (limit) aggregate process in $\mathcal{L}^2$ but not in $\mathcal{L}^1$. In so doing, we emphasize that complex analysis is required to obtain the convergence of the series $\sum_{k \geq 1}a_k$. Then we provide a characterization of the (limit) aggregate long-memory process in the presence of individual random AR(1) processes.

\subsection{The function $m$}

\noindent Before showing the convergence of the series of autoregressive parameters, we need an intermediate result regarding the $m$ function, and especially to make clear the relationship between the sum of the moving average coefficients $\sum_{k \geq 1}u_k$ and $m(r)$ when $r \rightarrow 1^{-}$. Proposition \ref{Proposition Equiv LimSum} clarifies this link and turns to be extremely useful when characterizing the infinite sum of the autoregressive coefficients of a (limit) aggregate long-memory process.\footnote{Note obviously that $\sumkzero u_k  = + \infty$ if $\lim_{\ r\rightarrow 1^{-}}m(r)=+\infty$.}

\begin{proposition}
\label{Proposition Equiv LimSum}
$\sumkone u_k = +\infty$
if and only if\/ $\lim_{\ r\rightarrow 1^{-}}m(r)=+\infty$
\end{proposition}

\noindent Proof : See Appendix A.


\subsection{The series $a(z) = \displaystyle\frac{m(z)}{1+m(z)}$ and the convergence of $\sum_{k \geq 1}a_k$}

\noindent One key issue to study the convergence of the autoregressive coefficients is that the existence of the limit, $\lim_{r\rightarrow 1^{-}}{a(r)}$,
does not insure its equality to $\sum_{k\geq 1} a_k$ or even the convergence of this sum.\footnote{For instance, consider the power series defined for $|r|<1$, $h(r) = \underset{k=0}{\overset{\infty}{\sum}}(-1)^kr^k$. Then $h(r) = \frac{1}{1+r}$ and $\underset{r \rightarrow 1^{-}}{\lim} h(r) = \frac{1}{2}$. However $\underset{0}{\overset{\infty}{\sum}}(-1)^k$ is not convergent.}
Supplementary {\it Tauberian\/} conditions are needed for this kind of results (Hardy, 1949; Titchmarsh, 1939; Korevaar, 2004). In particular, the convergence of this series has to be studied by making use of complex analysis, especially in the case of a (limit) aggregate long-memory process.
\\ \\
\noindent All of the proofs related to Theorem 1 are gathered in Appendix B. To summarize, the proof of the convergence of the autoregressive coefficients proceeds as follows. We first define $m$ as an analytic function of the complex variable $z$ within the open disc around 0 of radius one, $D(0,1)$, and rewrite $m(z)$ as an integral on $[0,1)$ with parameter $z \in D(0,1)$. This allows showing that it can be continuously extended to $\overline{D(0,1)}\backslash \{1\}$  (see Lemma 1 in Appendix B).  Then a second lemma (see Lemma 2 in Appendix B) proves that the function $a(z) = \frac{m(z)}{1+m(z)}$ is well defined in the disc $D(0,1)$ (i.e., the denominator does not vanish) and then it can be extended to the whole closed disc $\overline{D(0,1)}$.
For this purpose, two cases must be studied, according to the nature of the series $\sum u_k$: if $\sum u_k<\infty$ the function $m$ is continuous in the closed disc, so is $a$ ; and if
$\sum u_k = \infty$ then $|m(z)|\rightarrow \infty$ when $z\rightarrow 1$ in $D(0,1)$, therefore $a(z)\rightarrow 1$ when $z\rightarrow 1$ in $D(0,1)$.  Finally, a third lemma shows that the analytic function $a$ is univalent and provides the use of a Tauberian condition (see Lemma 3 in Appendix B). This then allows to prove the convergence of the series $\sum a_{n}$ (Theorem 1).

\begin{theorem}
Let $\left\{a_k, k \geq 1\right\}$ denote the sequence defined in Proposition \ref{Proposition Aggr Dyn}. The series $\sumkone a_{k}$ is convergent and 
\begin{equation*}
\sum_{k=1}^{+\infty} a_{k}=\lim_{r\rightarrow 1-}\frac{m(r)}{1+m(r)}
\end{equation*}
\label{Theorem_main}
\end{theorem}

\noindent Proof : See Appendix B.
\newline \newline
\noindent Taking Proposition \ref{Proposition Equiv LimSum} and Theorem \ref{Theorem_main}, it is then straightforward to determine the sum of the autoregressive coefficients.

\begin{proposition}
\label{proposition_sum_autoregressive_coefficients}
Let $\left\{a_k, k \geq 1\right\}$ denote the sequence defined in Proposition \ref{Proposition Aggr Dyn}.
The sum of the autoregressive coefficients, $\sumkone a_{k}$, equals one if and only if $\lim_{r\rightarrow 1^{-}}{m(r)}=+\infty$ or equivalently if and only if $\sumkone u_k = +\infty$.
\end{proposition}

\noindent Three points are worth commenting. First, if the limit aggregate process is a second-order stationary process (Proposition \ref{Proposition Aggr Dyn}) and the series of its moving average coefficients is absolutely summable, then the sum of the autoregressive coefficients is less than one. Notably, this result obtains with classical time series results whereas there is a need of complex analysis when the series of the moving average coefficients is not absolutely summable and one studies what happens at the pole $z = 1$. A second and related point is that the behavior of the series of the autoregressive coefficients depends on whether the limit second-order stationary aggregate process belongs to $\mathcal{L}^1$ or not. Consequently, as explained below, this provides a new characterization of a (limit) aggregate long-memory process as a result of the aggregation of random AR(1) processes. Third, as stated in Corollary \ref{corollary_zero}, the function $\Phi(z) = 1-a(z)$ admits only one zero on $\overline{D(0,1)}$.

\begin{corollary}
\label{corollary_zero}
Let $\Phi(z) = 1-a(z)$. Then $z=1$ is a zero of the function $\Phi$, which is defined on $\overline{D(0,1)}$,  if and only if $\sumkone u_k = +\infty$.
\end{corollary}

Proof : This is a straightforward implication of Theorem \ref{Theorem_main}.
\\
\\
\noindent Corollary \ref{corollary_zero} establishes that the representation of the stochastic process defined by the limit of the aggregate process $(X_{N,t})_{t \in \mathbb{Z}}$ admits a unit root whereas the aggregate process is weakly stationary. There is a one-to-one relationship between the long memory property and the presence of a unit root of the infinite autoregressive limit aggregate process.

\subsection{Time series implications}

\noindent Taking Theorem 1 and Proposition \ref{proposition_sum_autoregressive_coefficients}, we are now in a position to provide a new characterization of a (limit) aggregate long-memory process that results from the aggregation of individual random AR(1) processes.

\begin{theorem}
\label{long_memory_process}
Suppose that Assumptions 1 and 2 hold true. Let $(X_{t}, t \in \mathbb{Z})$ denote a long-memory process with the following Wold decomposition
\begin{eqnarray*}
X_{t} &=&\sum_{k=0}^{\infty }u_{k}\ \epsilon_{t-k}
\end{eqnarray*}
where the $u_k$ terms are nonnegative and $\underset{k=0}{\overset{\infty}{\sum}}u_k = \infty$. Then the $a_j$ terms of the equivalent infinite autoregressive representation,
\begin{eqnarray*}
X_{t} &=&\sum_{k=1}^{\infty }a_{k} X_{t-k}+\epsilon_{t},
\end{eqnarray*}
which are defined from $a_{k+1} = u_{k+1} -\sum_{r=1}^{k}a_{r}u_{k-r+1}$ and $a_1 = u_1$, satisfy $\sum_{k=1}^\infty a_k = 1$.
\end{theorem}

\noindent Several points are worth discussing. First, the (limit) aggregate long-memory process does not belong to the class of ARFIMA or self-similar stochastic processes.\footnote{A continuous-time stochastic process $(Y_t)$ is said to be self-similar with self-similarity paramter $H$, if for any sequence of time points $t_1$,$\cdots$,$t_k$ and any positive constant $a$, $c^{-H}(Y_{at_1},\cdots,Y_{at_k})$ has the same distribution as $(Y_{t_1},\cdots,Y_{t_k})$.} In particular, it is not possible to express the moving average coefficients of a fractional ARIMA process (by binomial expansion) such that they match the moving average weights of the long-memory process defined in Theorem \ref{long_memory_process}. In the same respect, matching the aggregate moving average coefficients of the (limit) aggregate process requires non-constant values of the self-similarity parameter $H$ (Beran, 1994; Beran et al., 2013). Second, the fact that the sum of the autoregressive coefficients of the (limit) aggregate process is equal to one is also consistent with the standard definition of long-memory processes (Beran, 1994; Beran et al., 2013), i.e. there exists a real number $\alpha \in (0,1)$ and a constant $c_{\rho}$ such that $\underset{k \rightarrow \infty}{\lim} \frac{\rho(k)}{c_{\rho}k^\alpha} = 1$ where $\rho(k) = \gamma_X(h)/\gamma_X(0)$  (with $\gamma_X(h) = \ex\left[\varphi^h/(1-\varphi^2)\right]$) is the autocorrelation of order $k$ of the (limit) aggregate stochastic process $(X_t)$.
\\ \\
\noindent Third, as pointed out by Beran (1994), observing long-range dependence in an aggregate time series (e.g., at the macro-level) does not necessarily mean that this is due to the genuine occurrence of long memory in the individual series (micro-level). This might be induced artificially by aggregation. Say differently, identifying the source of long memory would require to look carefully at the behavior of the possibly unobservable individual series. Fourth, Theorem \ref{long_memory_process} has some implications, which are beyond the scope of this paper, regarding some aggregation problems often encountered in the theoretical and empirical aggregation research (Pesaran and Chudik, 2014), namely the derivation of the macro dynamics from heterogenous individual dynamics and the identification and estimation of certain distributional features of the micro parameters from aggregate relations (disaggregation problem).\footnote{See Lewbel (1994), and Jondeau and Pelgrin ( 2014b).} For instance, in the case of the disaggregation problem when individual series are not available, since the autoregressive specification in Proposition \ref{Proposition Aggr Dyn} displays an infinite parametrization and the econometrician has only finitely many observations, one might proceed with a finite-parameter approximation for the infinite lag distribution (Sims, 1971, 1972; Faust and Lepper, 1997; P\"{o}tscher, 2002) and might account possibly for the constraint $a(1) =1$. Notably, finite-parameter approximations might lead to near observational aggregate processes with different statistical properties and thus to the incorrect identification of distributional features of the micro parameters.\footnote{The convergence of the estimates of the approximation is not sufficient to guarantee the convergence of some functions of those parameters---pointwise convergence does not imply (locally) uniform convergence.} Other applications might concern, among others, the estimation of the aggregate dynamics in the presence of unobserved heterogeneity when individual series are not available, the derivation of the macrodynamics in heterogeneous agents models, or the reconciliation of micro (panel) evidence and macro facts (e.g., the purchasing power parity, the nominal price rigidity).
\\ \\
\noindent Finally it is worth noting that Theorem \ref{long_memory_process} applies in a broader context than the aggregation of random AR(1) stochastic processes. Indeed, any stochastic process defined (almost surely) in $\mathcal{L}^2$ that admits a Wold decomposition with decreasing and nonnegative moment moving average coefficients and does not belong to $\mathcal{L}^1$ displays a unit root. This comes from the formal identity between the generating functions of $(a_k)$ and $(u_k)$, $\sum_{k\geq 1}a_{k}z^{k}
=\sum_{k\geq 1}u_{k}z^{k}-\sum_{k\geq 1}a_{k}z^{k}\sum_{k\geq 1}u_{k}z^{k}$. 

\section{Examples}

In this section, we review three examples, namely the Beta $B(p,q)$ distribution (of the first kind), the uniform distribution with $p = q = 1$, and the random polynomial density.

\paragraph{Beta distribution}
Following Morris and Granger (1976),  Granger (1980) and Gon\c{c}alves and Gouri\'{e}roux (1994), we assume that $\varphi$ is Beta-distributed.\footnote{Note that Gon\c{c}alves and Gouri\'{e}roux (1994) study extensively the aggregation of individual AR(1) processes in which $\varphi$ is Beta-distributed (after an homothety) and provide a discussion regarding aggregate long-memory processes. However they do not consider the sum of the autoregressive coefficients.}
\begin{equation*}
B(p,q ; x) = \frac{\Gamma(p+q)}{\Gamma(p)\Gamma(q)}
       x^{p-1}(1-x)^{q-1}\un_{[0,1)}(x),\ p > 0,\ q > 0.
\end{equation*}
In this respect, to the best of our knowledge, Proposition \ref{Beta_distribution} provides a new characterization of the series of autoregressive coefficients (i.e., the persistence of the (limit) aggregate process $(X_t)$).
\begin{proposition}
\label{Beta_distribution}
Suppose that $\varphi$ is Beta-distributed (of the first kind). Given the disaggregate model defined in Eqs. (\ref{xi,t})--(\ref{vi,t}), the series of the autoregressive coefficients of the limit aggregate process defined in Proposition \ref{Proposition Aggr Dyn} is given by:
\begin{itemize}
\item[-] If $q > 1$, then $\underset{k=1}{\overset{\infty}{\sum}} a_k = \frac{p}{p+q-1}$;
\item[-] If $q \leq 1$, then $\underset{k=1}{\overset{\infty}{\sum}} a_k = 1.$
\end{itemize}
\end{proposition}

\noindent Proposition \ref{Beta_distribution} can be shown as follows. 
Taking the integral form of $m(r)$ we have
\begin{equation*}
1+m(r)=\frac{\Gamma(p+q)}{\Gamma(p)\Gamma(q)}
            \int_0^1\frac{x^{p-1}(1-x)^{q-1}}{1-rx}dx,
\end{equation*}
and thus (by the monotone convergence theorem)
\begin{equation*}
\lim_{\ r\rightarrow 1^{-}}(1+m(r)) = \frac{\Gamma(p+q)}{\Gamma(p)\Gamma(q)}
            \int_0^1{x^{p-1}}{(1-x)^{q-2}}dx.
\end{equation*}
Consequently, this is a convergent integral if and only if $q >1$. In this case,
\begin{equation*}
\lim_{\ r\rightarrow 1^{-}}(1+m(r)) = \frac{\Gamma(p+q)}{\Gamma(p)\Gamma(q)}
            \frac{\Gamma(p)\Gamma(q-1)}{\Gamma(p+q-1)}=1+\frac{p}{q-1}
           \end{equation*}
and
$$\sum_{k=1}^{+\infty} a_{k}=\lim_{r\rightarrow 1-}\frac{m(r)}{1+m(r)}=\frac{p}{p+q-1}.$$
\noindent On the other hand, it follows that  $\lim_{\ r\rightarrow 1^{-}}(1+m(r)) =+\infty$ if and only if $q \leq 1$, and $\sum_{k=1}^{+\infty} a_{k}=1$ (Theorem 1).

\paragraph{Uniform distribution}

\noindent We now assume that $p = q =1$ such that the random variable $\varphi$ is uniformly distributed over the interval $[0,1)$ and non central moments are given by
$u_k=\frac{1}{k+1}$.

\begin{proposition}
Suppose that $\varphi$ is uniformly distributed over the interval $[0;1)$. Given the disaggregate model defined in Eqs. (\ref{xi,t})--(\ref{vi,t}), the autoregressive coefficients of the limit aggregate process defined in Proposition \ref{Proposition Aggr Dyn} are given by:
\begin{eqnarray*}
a_k = \frac{|I_k|}{k!}
\end{eqnarray*}
where $I_k = \int_0^1x(x-1)\cdots(x-k+1)dx$ has the same sign as $(-1)^{k-1}$ for $k\geq 1$. Moreover, $\sum_{k=1}^{\infty} a_k = 1$.
\label{uniform_distribution}
\end{proposition}

\noindent Proposition \ref{uniform_distribution} can be shown by using either Theorem 1 or a new lemma provided in Appendix C. Notably, the coefficients of the series $a(z)=m(z)/(1+m(z))$ can be computed as follows.
First, the generating moment series is
\begin{equation*}
1+m(z) = \sum_{k \geq 0}u_k z^k = \sum_{k\geq 0}\frac{z^k}{k+1} = - \frac{\log(1-z)}{z}
\end{equation*}
where $|z|<1$ so that
\begin{equation*}
a(z) = \frac{m(z)}{1+m(z)} = 1+\frac{z}{\log(1-z)}.
\end{equation*}
Second the expression of $\frac{z}{\log(1-z)}$ is derived by a power series development of the function $a$. Indeed, using the notation $\psi(z)=\log(1-z)$, one has
\begin{equation*}
\int_0^1 e^{x\psi(z)}dx = \frac{e^{\psi(z)}-1}{\psi(z)} = -\frac{z}{\log(1-z)}
\end{equation*}
where the power series development of $e^{x\psi(z)} = (1-z)^x$  is defined to be:
\begin{equation*}
(1-z)^x = \sum_{n\geq 0} {x\choose n} (-1)^nz^n
    = 1+\sum_{n\geq 1}\frac{(-1)^n}{n!}z^nx(x-1)\cdots(x-n+1).
\end{equation*}
Since this series converges absolutely for $|z|<1$ and uniformly for $x$ in $[0,1]$, one obtains
\begin{equation*}
\int_0^1 e^{x\psi(z)}dx  = \int_0^1 (1-z)^xdx  = 1+\sum_{n\geq 1}\frac{(-1)^n}{n!}z^nI_n
\end{equation*}
where $I_n = \int_0^1x(x-1)\cdots(x-n+1)dx$ for $n\geq 1$. Finally,
\begin{equation*}
\int_0^1 e^{x\psi(z)}dx  = -\frac{z}{\log(1-z)} = 1 - \sum_{k\geq 1}\frac{1}{k!}z^k|I_k|
\end{equation*}
and thus
\begin{equation*}
a(z) =  \sum_{k\geq 1}\frac{|I_k|}{k!}z^k.
\end{equation*}
Several points are worth commenting. On the one hand, the non-negativeness of the autoregressive coefficients for the uniform distribution allows for the use of a standard Hardy-Littlewood Tauberian result and thus might not require Theorem 1. However the non-negativeness of the autoregressive coefficients is not proved for all probability distributions $\mu$ with support $[0,1)$. On the other hand, using Theorem 1, since the moving average coefficients of the (limit) aggregate process constitute a harmonic series, the corresponding series diverges and thus the sum of the autoregressive coefficients (of the limit aggregate process) equals one. Finally, Proposition \ref{uniform_distribution} extends the result of Linden (1999) in which the behavior of the aggregate process is studied with the autocorrelation function.

\paragraph{Random polynomial density}


\noindent We consider as a last example the case of the polynomial aggregated AR(1) model. More specifically, we suppose that
$\varphi$ has a continuous distribution over $[0,1]$ that can be represented by a polynomial of degree $d\geq 1$ :
\begin{eqnarray*}
f(\varphi) = \sum_{s=0}^d c_s \varphi^s\un_{[0,1]}(\varphi)
\end{eqnarray*}
where $\sum_{s=0}^{d} \frac{c_s}{s+1} = 1$ (the density has to be integrated to one) and $f$ is non-negative in $[0,1]$. One key feature of the polynomial density function relative to the Beta distribution is that it can be multi-modal (for polynomial density of third order or above). Similarly to the Beta distribution, two cases are considered in Proposition \ref{random_polynomial_distribution}: the sum of the autoregressive coefficients does depend on whether $f(1) = 0$ or $f(1) > 0$ with $f(1)=c_0+c_1+\cdots +c_d$ is the value of the density at $x=1$.

\begin{proposition}
\label{random_polynomial_distribution}
Suppose that $\varphi$ has a polynomial density function of order $d$. Given the disaggregate model defined in Eqs. (\ref{xi,t})--(\ref{vi,t}), the series of the autoregressive coefficients of the limit aggregate process defined in Proposition \ref{Proposition Aggr Dyn} is given by:
\begin{itemize}
\item[-] If $f(1) = 0$, then
\begin{eqnarray*}
a(1) = \sum_{k=1}^\infty a_k = 1 -  \frac{1}{\underset{n=0}{\overset{d-1}{\sum}}\underset{k=0}{\overset{n}{\sum}} \frac{c_{n-k}}{n+1}};
\end{eqnarray*}
\item[-] If $f(1) > 0$, then $\underset{k=1}{\overset{\infty}{\sum}} a_k = 1.$
\end{itemize}
\end{proposition}

\noindent Indeed, starting from the polynomial density function of order $d$, the non-central moments $u_k = \ex\left[\varphi^k\right]$, $k\geq 0$,  are given by:
\begin{eqnarray*}
\ex\left[\varphi^k\right] = \sum_{s=0}^d \frac{c_s}{s+k+1},
\end{eqnarray*}
and thus the generating moment series,
\begin{equation*}
1+m(z) = \sum_{k\geq 0}u_k z^k = \sum_{k=0}^{\infty}\Bigl(\sum_{s=0}^d \frac{c_s}{s+k+1}\Bigr)z^k,
\end{equation*}
is convergent at least for $|z|<1$. Therefore, one needs to study the convergence of the series
\begin{equation*}
1+m(1)=\sum_{k=0}^{\infty}\sum_{s=0}^d \frac{c_s}{s+k+1}
          = \lim_{K\rightarrow \infty}\sum_{k=0}^{K}\sum_{s=0}^d \frac{c_s}{s+k+1}
\end{equation*}

\noindent The terms of the double sum $S(K) = \sum_{k=0}^{K}\sum_{s=0}^d \frac{c_s}{s+k+1}$ form an array with
$K+1$ rows and $d+1$ columns. Without loss of generality, we can suppose $K>d$. Let $n=s+k$, we have $n-d\leq k\leq n$.
Taking diagonal sums along the lines $s+k=0,1,\ldots,d+K$ of the array, one can write
$S(K)=S_1(K)+S_2(K)+S_3(K)$ with
\begin{eqnarray*}
S_1(K) &=& \sum_{n=0}^{d}\sum_{k=0}^{n} \frac{c_{n-k}}{n+1}
       = c_0+\frac{c_0+c_1}{2}+\cdots +\frac{c_0+c_1+\cdots +c_{d}}{d+1} \\
S_2(K) &=& \sum_{n=d+1}^{K}\sum_{k=n-d}^{n} \frac{c_{n-k}}{n+1}
       = \sum_{n=d+1}^{K}\frac{1}{n+1}\sum_{k=n-d}^{n} c_{n-k}
       = f(1)\sum_{n=d+1}^{K}\frac{1}{n+1}
\end{eqnarray*}
and
\begin{eqnarray*}
S_3(K) &=& \sum_{n=K+1}^{d+K}\sum_{k=n-d}^{n} \frac{c_{n-k}}{n+1}
            = \sum_{n=K+1}^{d+K}\frac{1}{n+1}\sum_{k=n-d}^{n}c_{n-k} \\
       &=& \frac{c_1+c_2+\cdots +c_{d}}{K+2}+\frac{c_2+\cdots +c_{d}}{K+3}+\cdots +\frac{c_d}{K+D+1}.
\end{eqnarray*}
The first sum, $S_1(K)$, is finite and independent of $K$. The second sum, $S_2(K)$, is clearly convergent as $K \rightarrow \infty$ if and only if $f(1) = 0$. Finally, the third sum, $S_3(K)$, is finite and its limit is zero when $K$ approaches infinity. In this respect, the generating moment series $1+m(1)$ converges if and only if $f(1)=0$. In this case its limit is $S_1=S_1(K)$ with
$$
S_1(K)=\displaystyle c_0+\frac{c_0+c_1}{2}+\cdots +\frac{c_0+c_1+\cdots +c_{d}}{d+1}
$$
In contrast, when $f(1)>0$ the series $\sum_{k\geq 0}u_k$ diverges and thus the sum of the autoregressive coefficients equals one. Consequently, the (limit) aggregate process displays long-range dependance.
\\ \\
\noindent In this respect, Proposition \ref{random_polynomial_distribution} provides formally the sufficient condition discussed by Chong (2006), i.e. $f(1) > 0$ is sufficient to establish the long memory properties of the (limit) aggregate process. Moreover, if we assume that the polynomial density function is of order zero and $c_0 = 1$, we end up with the uniform distribution (with $f(1) > 0$) and thus Proposition \ref{uniform_distribution}. On the other hand, if we assume that the polynomial density function is of the form $ax(1-x)$ for $a\neq 0$ and $x \in [0,1)$, then $f(1) = 0$ and the long-run persistence is given by $a(1)$. In contrast, if the probability density function has a support on $[0;1]$, the aggregation process leads to a generalized integrated process (Lin, 1991; Granger and Ding, 1996). As a final remark, it is worth emphasizing that any distribution such that the generating function of the sequence of the moments is not convergent leads to a long-memory process characterized by Theorem 2.

\section{Conclusion}

\noindent In this paper, we study the aggregation of individual random AR(1) processes under the assumptions that some common factors exist, the (positive) random autoregressive coefficient takes values in $[0,1)$ with probability distribution $\mu$, and all of the noncentral moments exist. Notably we show by making use of complex analysis that sum of the autoregressive coefficients equals one in the presence of limit aggregate long memory processes: the divergence of the series of noncentral moments is fully equivalent to the presence of a unit root in the autoregressive representation of the stochastic process defined by the limit of the aggregate process. We then illustrate our results using some prominent examples of distribution for the aggregation of random autoregressive AR(1) processes. This provides some new insights that might deserve some empirical applications and some theoretical developments, as for instance the disaggregation problem.

\newpage

\section*{Appendix: Proofs}

\subsection*{Appendix A}

\noindent Proof of Proposition \ref{Proposition Equiv LimSum}: Taking that the sequence of the moments
$$u_{k} = \ex(\varphi^k) = \int_{[0,1)}x^kd\mu(x)$$  is positive and decreasing, one has $\lim_{k \rightarrow +\infty}u_k=0$ (by monotone convergence applied to the sequence $(x^k)_{k\geq 1}$ and thus the radius of convergence of the series $\sum_{k\geq 1}u_{k}z^{k}$ is  at least 1.

\noindent Let $(r_n)$ an increasing sequence in $[0,1)$ such that $\lim_{n\rightarrow +\infty} r_n=1$.
For $x\in [0,1)$ the functions $f_{n}(x) = \displaystyle\frac{r_{n}x}{1-r_{n}x}$
are positive therefore by Fatou's Lemma we get
\begin{equation*}
\lim \inf_{n}m(r_n)=\lim \inf_{n}\int_{[0,1)}f_{n}d\mu(x)\geq \int_{[0,1)}\lim \inf_{n}f_{n}d\mu(x)=\int_{[0,1)}\frac{x}{1-x}d\mu(x)=
\sum_{k=0}^{+\infty }u_{k}
\end{equation*}

\noindent Thus if $\sum u_{k}$ is divergent then $m(r)\rightarrow+\infty$ when $r\rightarrow 1^{-}$.
 Moreover by the Abel Theorem (Titchmarsh, 1939, pp.9-10) if $\sum u_{k}$ is convergent then $\lim_{r\rightarrow 1^{-}} m(r)$ is finite and
$\lim_{r\rightarrow 1^{-}} m(r)=\sum u_{k}$. The equivalence then follows.

\subsection*{Appendix B}

In Appendix B, we provide the proof of Theorem 1. In so doing, we proceed with three lemmas. Notably we first define $m$ as an analytic function of the complex variable $z$ in the disc $D(0,1)$ and rewrite $m(z)$  as an integral on $[0,1)$ with parameter $z \in D(0,1)$ (Lemma 1). We then show that it can be continuously extended to $\overline{D(0,1)}\backslash \{1\}$.  Then Lemma 2 proves that
the function $a(z) = \frac{m(z)}{1+m(z)}$ is well defined in the disc $D(0,1)$ (i.e., the denominator does not vanish) so that it can be extended to the whole closed disc $\overline{D(0,1)}$.
In this respect, two cases must be studied, according to the nature of the series $\sum u_k$: if $\sum u_k<\infty$ the function $m$ is continuous in the closed disc, so is $a$ ; and if
$\sum u_k = \infty$ then $|m(z)|\rightarrow \infty$ when $z\rightarrow 1$ in $D(0,1)$, therefore $a(z)\rightarrow 1$ when $z\rightarrow 1$ in $D(0,1)$.  Finally, Lemma 3 proves that the analytic function $a$ is shown to be univalent and provides the use of a Tauberian condition. Then Theorem 1 is proven.

\setcounter{lemma}{0}

Let $m$ the function defined in the open disc $D(0,1)=\{z\in\mathbb{C} \text{ with } \vert z \vert <1\}$ by
\begin{equation*}
m(z)=\sum_{n=1}^{+\infty }u_{n}z^{n}
\end{equation*}
This function is analytic.  Moreover,  by a classical theorem,  on the boundary
$C(0,1)=\{z\in\mathbb{C} \text{ with } \vert z \vert =1\}$  the series
$\sum_{n\geq 1}u_{n}e^{int}$ is convergent  for all  $t\in ]0,2\pi[$ and by Abel's theorem :
\begin{equation*}
m(e^{it})=\sum_{n=1}^{+\infty }u_{n}e^{int}=\lim_{r\rightarrow 1-}\sum_{n=1}^{+\infty }u_{n}r^{n}e^{int}
\end{equation*}
Thus the function $m$ is defined in $\overline{D(0,1)}\backslash \{1\}=\{z\in\mathbb{C} \text{ with } \vert z \vert \leq1,\ z\neq1\}$.
By  positivity we  have
\begin{eqnarray*}
\sum_{n=1}^{+\infty }E(\varphi^{n})r^{n}=E(\sum_{n=1}^{+\infty }\varphi^{n}r^{n})
\end{eqnarray*}
for all $0\leq r<1$. Thus for  $r\in[0,1)$ we get
\begin{equation*}
m(r)=\sum_{n=1}^{+\infty }u_{n}r^{n}=E(\frac{\varphi r}{1-\varphi r})=\int_{[0,1)}\frac{rx}{1-rx}d\mu(x)
\end{equation*}

\begin{lemma}
For $z\in \overline{D(0,1)}\backslash \{1\}$
\begin{equation*}
m(z)=\int_{[0,1)}\frac{zx}{1-zx}d\mu(x)
\end{equation*}
\label{representation_m}
\end{lemma}

Proof: First we prove that the function $z\mapsto \int_{[0,1)}\frac{zx}{1-zx}d\mu(x)$ is analytic in the open disc $D(0,1)$,  continuous in $\overline{D(0,1)}\backslash \{1\}$ and thus is defined.
\newline \newline
\noindent Let $K$  denote a compact set in $\overline{D(0,1)}\backslash \{1\}$. This compact is
included in $\overline{D(0,1)}\backslash D(1,\varepsilon )$ with
$\varepsilon >0$ so if $z\in K$ we get $zx\in \overline{D(0,1)}\backslash
D(1,\varepsilon )$ for all $x\in [0,1)$, and thus
\begin{equation*}
\left| 1-zx\right| \geq \varepsilon.
\end{equation*}
\noindent Therefore for all $z\in K$ and $x\in ]0,1[$
\begin{equation*}
\left| \frac{zx}{1-zx}\right| \leq \frac{1}{\varepsilon }x.
\end{equation*}
\noindent This proves the continuity of the function $z\mapsto \int_{[0,1)}\frac{zx}{1-zx}d\mu(x)$
over $K$. Moreover, by analyticity of $z\mapsto \frac{zx}{1-zx}$ for all $x\in
[0,1)$ and the previous boundness condition, the function $z\mapsto \int_{[0,1)}\frac{zx}{1-zx}d\mu(x)$
is analytic in $D(0,1)$.

Finally we note that this function coincides with $m(r)$, $r\in [0,1[$, thus by analytic continuation, we obtain
\begin{equation*}
m(z)=\int_{[0,1)}\frac{zx}{1-zx}d\mu(x)
\end{equation*}
for all  $z\in D(0,1)$, and for $z=e^{it}$ with $t\neq 0$
\begin{equation*}
m(e^{it})=\lim_{r\rightarrow 1-}\sum_{n=1}^{+\infty }u_{n}r^{n}e^{int}
=\lim_{r\rightarrow 1-} \int_{[0,1)}\frac{re^{it}x}{1-re^{it}x}d\mu(x)
= \int_{[0,1)}\frac{e^{it}x}{1-e^{it}x}d\mu(x)
\end{equation*}
$\square$

\textbf{Extension by continuity of the function $\frac{m}{1+m}$}
\begin{lemma}
The function $\frac{m}{1+m}$
can be extended to a continuous function over $\overline{D(0,1)}$.
\label{Lemma_continuous}
\end{lemma}

Proof: The function $1+m$ doesn't vanish in $\overline{D(0,1)}$ because
\begin{equation*}
1+m(z)=\int_{[0,1)}\frac{1}{1-zx}d\mu(x)
\end{equation*}

and if $z=a+ib$ with  $a^{2}+b^{2}\leq1$ then $ax<1$ for all $x\in [0,1)$ so
\begin{equation*}
\text{Re}(1+m(z))=\int_{[0,1)}\frac{1-ax}{(1-ax)^{2}+b^{2}x^{2}}d\mu(x)>0
\end{equation*}

The function $\frac{m}{1+m}$ is therefore defined and continuous in
$\overline{D(0,1)}\backslash \{1\}$ and analytic in $D(0,1)$.
\par\smallskip

In order to study the continuity at the point $z=1$, we need to consider two cases :
\begin{enumerate}
\item If $\sum u_{k}$ is convergent then the series $\sum u_{k}z^{k}$ is
normally convergent in $\overline{D(0,1)}$ and the function $m$ is
continuous in  $\overline{D(0,1)}$ and consequently so is $\frac{m}{1+m}$.
\item If $\sum u_{k}$ is divergent it is sufficient to prove that
$|m(z)|\rightarrow +\infty $ when $z\rightarrow 1$. Indeed
\begin{equation*}
\frac{m(z)}{1+m(z)}=\frac{1}{1+\frac{1}{m(z)}}\rightarrow 1
\end{equation*}
and we extend $\frac{m}{1+m}$ by $1$ at the point $z=1$.
\end{enumerate}
\noindent We  now show that $|m(z)|\rightarrow +\infty $ when $z\rightarrow 1.$ In so doing, consider a sequence of points $z_{k}=a_{k}+ib_{k}$ with
$a_{k}\rightarrow 1,b_{k}\rightarrow 0$ and $a_{k}^{2}+b_{k}^{2}\leq 1$.
Then
\begin{equation*}
\text{Re(}m(z_{k}))=-1+\int_{[0,1)}\frac{1-a_{k}x}{(1-a_{k}x)^{2}+b_{k}^{2}x^{2}}d\mu(x).
\end{equation*}
\noindent As $a_{k}x\leq 1$ for all integer $k$ and all $x\in [0,1)$ the functions
\begin{equation*}
f_{k}:x\mapsto \frac{1-a_{k}x}{(1-a_{k}x)^{2}+b_{k}^{2}x^{2}}
\end{equation*}
are positive in $[0,1)$. Using Fatou's Lemma, we get
\begin{equation*}
\lim \inf_{k}\int_{[0,1)}f_{k}d\mu(x)\geq \int_{[0,1)}\lim \inf_{k}f_{k}d\mu(x)
\end{equation*}
where
\begin{equation*}
\lim \inf_{k}f_{k}(x)=\lim f_{k}(x)=\frac{1}{1-x}.
\end{equation*}
It is therefore sufficient to remark that if $\sum u_{k}$ is divergent then
\begin{equation*}
\int_{[0,1)}\frac{1}{1-x}d\mu(x)\geq
\int_{[0,1)}\frac{x}{1-x}d\mu(x)=
\sum_{n=0}^{+\infty }u_{n}=
+\infty.
\end{equation*}
\noindent So for all sequence $(z_{k})$ converging to 1 in the disc, we have
$Re(m(z_{k}))\rightarrow +\infty$ and thus
$|m(z_{k})|\rightarrow +\infty $.
$\square$
\\ \\
\noindent \textbf{Convergence of an analytic function $D(0,1)$, which is continuous on $\overline{D(0,1)}$}
\\ \\
\noindent \textbf{Lemma \ref{Lemma_convergence}} provides the use of a Tauberian condition, which turns to be crucial for the proof of Theorem 1.

\begin{lemma}
Let $f(z)=\sum_{k\geq 1}b_{k}z^{k}$ an analytic function on $D(0,1)$, continuous on $\overline{D(0,1)}.$
If $f$ is injective on $D(0,1)$ then the series
$\sum_{k\geq 1}b_{k}$ is convergent.
\label{Lemma_convergence}
\end{lemma}

\noindent Proof:  We proceed in two steps. First, we prove that
$$
\frac{1}{n}\sum_{k=1}^{n}k\left|b_{k}\right| \rightarrow 0.
$$
On the one hand, the function $f$ is analytic on $D(0,1)$ thus the image $U=f(D(0,1))$ is open
in $\mathbb{C}$ and $U$ is included in the compact set $f(\overline{D(0,1)})$ since $f$ is continuous on $\overline{D(0,1)}.$
\newline \newline
On the other hand, the function  $f$ being injective on $D(0,1)$ we have  $f'(z)\neq 0$ for all $z\in D(0,1)$ thus $f$ is a $C^{1}$ diffeomorphism between $D(0,1)$ and
$U$. By the change of variables formula we get
\begin{equation*}
\lambda (U)=\int_{D(0,1)}\left| f^{\prime }(x+iy)\right| ^{2}dxdy
\end{equation*}
and with the use of  polar coordinates we get the finite Lebesgue measure of  $U$ as the sum of the series
\begin{equation*}
\lambda (U)=\pi \sum_{n=1}^{+\infty}n|b_{n}|^{2}
\end{equation*}
The convergence of this last series now implies that $\frac{1}{n}\sum_{k=1}^{n}k\left|b_{k}\right| \rightarrow 0$.
\newline \newline
Finally, to verify this assertion let $N\geq 1$ and for $n>N$  write
$$
\frac{1}{n}\sum_{k=1}^{n}k\left| b_{k}\right|  = \frac{1}{n}%
\sum_{k=1}^{N}k\left|b_{k}\right| +\frac{1}{n}\sum_{k=N+1}^{n}k\left|b_{k}\right|.
$$
Then by Cauchy-Schwarz inequality we have
\begin{eqnarray*}
\frac{1}{n}\sum_{k=1}^{n}k\left| b_{k}\right| &\leq &\frac{1}{n}\sum_{k=1}^{N}k\left|b_{k}\right| +\frac{1}{n}%
(\sum_{k=N+1}^{n}k)^{1/2}(\sum_{k=N+1}^{n}k\left|b_{k}\right| ^{2})^{1/2} \\
&\leq &\frac{1}{n}\sum_{k=1}^{N}k\left|b_{k}\right| +\frac{1}{n}%
(\sum_{k=N+1}^{n}k\left|b_{k}\right| ^{2})^{1/2}
\end{eqnarray*}
and thus
\begin{equation*}
\lim \sup_{n\rightarrow +\infty }\frac{1}{n}\sum_{k=1}^{n}k
\left|b_{k}\right| \leq \frac{1}{n}(\sum_{k=N+1}^{n}k\left|b_{k}\right|
^{2})^{1/2}.
\end{equation*}
Since $\frac{1}{n}(\sum_{k=N+1}^{n}k\left|b_{k}\right|^{2})^{1/2}\rightarrow 0$
it follows that
$\frac{1}{n}\sum_{k=1}^{n}k\left|b_{k}\right| \rightarrow 0$.
\newline \newline
Taking this intermediate result, we are now in a position to prove {\bf Lemma \ref{Lemma_convergence}}. Indeed
let $t_{n}=\frac{1}{n}\sum_{k=1}^{n}kb_{k}$
(with $t_{0}=0$). We have
$
b_{n} =(t_{n}-t_{n-1})+\frac{1}{n}t_{n-1}.
$
The series $\sum_{n\geq 1}(t_{n}-t_{n-1})$ is convergent because
\begin{equation*}
\sum_{n=1}^{N}(t_{n}-t_{n-1})=t_{N}\rightarrow 0,
\end{equation*}
and thus the series $\sum_{n\geq 1}(t_{n}-t_{n-1})$ is Abel summable.
Since the series $\sum_{n\geq 1}b_{n}$ is also Abel summable by continuity of $f$, we get the Abel-summability of the series
$\sum_{n\geq 1}\frac{1}{n}t_{n-1}$.
Finally, since
\begin{equation*}
n\bigl(\frac{1}{n}t_{n-1}\bigr)=t_{n-1}\rightarrow 0,
\end{equation*}
it follows from the classical Tauber's theorem that the series $\sum_{n\geq 1}\frac{1}{n}t_{n-1}$ (and thus $\sum_{n\geq 1}b_{n}$) is convergent.
$\square$
\\ \\
\noindent \textbf{Proof of Theorem 1}
\\ \\
\noindent Using {\bf Lemma 2\/}, the function $f=a=\frac{m}{1+m}$ is defined and continuous on
$\overline{D(0,1)}$ and analytic in $D(0,1)$.
It follows from  {\bf Lemma 3\/} that it remains to prove that $f$ is
injective in $D(0,1)$. Taking
\begin{equation*}
\frac{m(z_{1})}{1+m(z_{1})}=\frac{m(z_{2})}{1+m(z_{2})}\Leftrightarrow
m(z_{1})=m(z_{2}),
\end{equation*}
it is sufficient to prove the injectivity of the function $m$ on $D(0,1)$.
Moreover it is sufficient to prove this injectivity on
$D(0,r) $ for all $0 \leq r < 1$ : if $m$ is not injective on $D(0,1)$
then there exists $z_{1}\neq z_{2}$ such that $m(z_{1})=m(z_{2})$ ; so $m$
is not injective on $D(0,r)$ where $r>\max (\left| z_{1}\right| ,\left|
z_{2}\right| )$.
\par\smallskip
For all $0 \leq r < 1$ the function $m$ is analytic on $D(0,r)$ and continuous
on $\overline{D(0,r)}$. Using Darboux's theorem (Burckel, 1979, p. 310), we could then establish the injectivity of
$m$ on $D(0,r)$ by showing that $m$ is injective on the circle of radius $r$ with center 0.
Indeed, let
\begin{equation*}
\varphi (t)={Re}(m(re^{it}))=
\int_{[0,1)}\frac{rx\cos (t)-r^{2}x^{2}}{1-2rx\cos (t)+r^{2}x^{2}}d\mu(x).
\end{equation*}
We have
\begin{equation*}
\varphi ^{\prime }(t)=\sin t\int_{[0,1)}\frac{-rx(1-r^{2}x^{2})}{\left(
1-2xr\cos t+r^{2}x^{2}\right) ^{2}}d\mu(x).
\end{equation*}
We see that the function $\varphi $ is decreasing on
$]0,\pi[$ et increasing on $]\pi ,2\pi[$.
It is symmetric across $\pi $ : we have
$\varphi (t)=\varphi(2\pi -t)$. Therefore the only points $t_{1}\neq t_{2}$ with $\varphi
(t_{1})=\varphi (t_{2})$ are the pairs $(t,2\pi -t)$ with $t\in [0,\pi[$.
Moreover we have
\begin{equation*}
{\text{Im}}(m(re^{it}))=\int_{[0,1)}\frac{xr\sin (t)}{1-2xr\cos
(t)+r^{2}x^{2}}d\mu(x)
\end{equation*}
hence
\begin{equation*}
{\text{Im}}(m(re^{it}))=-{\text{Im}}(m(re^{i(2\pi -t)}))
\end{equation*}
Since we have
\begin{equation*}
\sin (t)\int_{[0,1)}\frac{xr}{1-2xr\cos (t)+r^{2}x^{2}}d\mu(x)>0
\end{equation*}
for $t\in [0,\pi)$, we can't have
\begin{equation*}
{\text{Im}}(m(re^{it}))={\text{Im}}(m(re^{i(2\pi -t)}))
\end{equation*}
Therefore $t\mapsto m(re^{it})$ is injective.
By the Abel Theorem the sum of the series
$\sum a_{n}$  is equal to the limit
$\lim_{r\rightarrow 1}a(r)$
where $r\in [0,1)$
$\square $

\newpage

%
%
%
%
%
%

\end{document}